\documentclass[12pt]{amsart}
\usepackage{amssymb,latexsym, amsmath, amscd, array, graphicx,pb-diagram}

\swapnumbers \numberwithin{equation}{section}

\theoremstyle{plain}

\newtheorem{thm}{Theorem}[section]

\newtheorem{conjec}[thm]{Conjecture}
\newtheorem{prop}[thm]{Proposition}
\newtheorem{cor}[thm]{Corollary}

\theoremstyle{definition}
\newtheorem{defn}[thm]{Definition}

\newtheorem{question}[thm]{Question}

 \newcommand{\Wi}{\widetilde}

\def\C{{\mathbb C}}
\def\Z{{\mathbb Z}}

\def\R{{\mathbb R}}

\def\1{\hbox{\rm\rlap {1}\hskip.03in{\rom I}}}
\def\Bbbone{{\rm1\mathchoice{\kern-0.25em}{\kern-0.25em}
{\kern-0.2em}{\kern-0.2em}I}}

\long\def\forget#1\forgotten{} %

\newcommand\ver[1]{\marginpar{\tiny Changed in Ver \VER}}

\newcommand{\mc}{ \text {mc}}

\newcounter{notecounter}

\date{\today}

\begin{document}

\title[On Gromov's conjecture]{On Gromov's conjecture for  totally non-spin  manifolds}

\author[D.~Bolotov]{Dmitry Bolotov}

\author[A.~Dranishnikov]{Alexander  Dranishnikov$^{1}$} 
\thanks{$^{1}$Supported by NSF grant DMS-1304627}

\address{Dmitry Bolotov, Verkin Institute of Low Temperature Physics, Lenina Avenue, 47, Kharkov, 631103,
Ukraine}
\email{bolotov@univer.kharkov.ua}

\address{Alexander N. Dranishnikov, Department of Mathematics, University
of Florida, 358 Little Hall, Gainesville, FL 32611-8105, USA}
\email{dranish@math.ufl.edu}



\begin{abstract}
Gromov's conjecture states that for a closed $n$-manifold $M$ with positive scalar curvature the macroscopic dimension of its universal covering $\Wi M$ satisfies the inequality $\dim_{mc}\Wi M\le n-2$~\cite{G2}.  
We prove that for totally non-spin $n$-manifolds the inequality $\dim_{mc}\Wi M\le n-1$ implies the inequality $\dim_{mc}\Wi M\le n-2$. This implication together with the main result of~\cite{Dr1} allows us to prove Gromov's conjecture
for totally non-spin $n$-manifolds whose fundamental group is a virtual duality group with $vcd\ne n$.

In the case of virtually abelian groups we reduce Gromov's conjecture for totally non-spin manifolds to the problem whether $H_n(T^n)^+\ne 0$. 
This problem can be further reduced to the $S^1$-stability conjecture for manifolds with  free abelian fundamental groups.
\end{abstract}

%%% ----------------------------------------------------------------------


%%% ----------------------------------------------------------------------
\maketitle
%%% ----------------------------------------------------------------------

\section {Introduction}
The notion of macroscopic dimension was introduced by M.
Gromov~\cite{G2} to study the topology of manifolds that admit a positive scalar curvature (PSC) metric.
We recall that the scalar curvature of a Riemannian $n$-manifold $M$ is a function
$ Sc_M:M\to\R$ which assigns to each point $x\in M$  two times the sum of the sectional curvatures over all 2-planes $e_i\wedge e_j$ in the tangent space
$T_xM$ at $x$ for some orthonormal basis $e_1,\dots,e_n$.
\begin{defn}
 A metric space $X$ has the macroscopic dimension $\dim_{\mc} X \leq k$ if
there is a uniformly cobounded proper map $f:X\to K$ to a $k$-dimensional simplicial complex.
Then $\dim_{mc}X=m$ where $m$ is minimal among $k$ with $\dim_{\mc} X \leq k$.
\end{defn}
\smallskip
We recall that a map of a metric space $f:X\to Y$ is uniformly cobounded if there is a uniform upper bound on the diameter of preimages $f^{-1}(y)$, $y\in Y$.

\

{\bf Gromov's Conjecture.} {\it The macroscopic dimension of the universal covering $\Wi M$ of  a closed  PSC  $n$-manifold $M$ satisfies the inequality $\dim_{mc}\Wi M\leq n-2$ for the metric on $\Wi M$ lifted from $M$}.

\

The main examples supporting Gromov's conjecture are $n$-manifolds of the form $M=N\times S^2$. They admit metrics with PSC in view of
the formula $Sc_{x_1,x_2}=Sc_{x_1}+Sc_{x_2}$ for the Cartesian product $(X_1\times X_2,\mathcal G_1\oplus\mathcal G_2)$
of two Riemannian manifolds $(X_1,\mathcal G_1)$ and $(X_2,\mathcal G_2)$ and the fact that while $Sc_N$ is bounded $Sc_{S^2}$ can be chosen
to be arbitrarily large. 
Note that the projection $p:\Wi M=\Wi N\times S^2\to\Wi N$ is a proper uniformly cobounded map to a $(n-2)$-dimensional manifold.
Hence, $\dim_{mc}\Wi M\le n-2$.

Since $\dim_{mc}X=0$ for every compact metric space, Gromov's conjecture holds trivially for simply connected manifolds. Thus, this conjecture  is about manifolds with nontrivial fundamental groups.
To what extend is Gromov's conjecture a conjecture about groups? This is the question that we are trying to answer.
We say that Gromov's conjecture holds for a group $\pi$ if it holds for manifolds with the fundamental group $\pi$. 
Thus, it makes sense to investigate Gromov's conjecture for classes of groups. Clearly, the conjecture holds true for all finite groups. This paper is an attempt to establish 
the conjecture for the class of virtual duality groups. 

Dealing with PSC manifolds one has to consider three different cases: the case of spin manifolds, almost spin manifolds, and totally non-spin
manifolds. We adopt the names {\em almost spin} for manifolds with a spin universal covering and {\em totally non-spin} for manifolds whose universal coverings are non-spin.

We note that in the case of spin manifolds (as well as almost spin) there is  index theory which provides a technique for attacking Gromov's conjecture. 
It was used to prove the conjecture for spin manifolds with a fundamental group satisfying the analytic Novikov conjecture and 
the Rosenberg-Stolz condition on injectivity of the real K-theory periodization map $per: ko_*(B\pi)\to KO_*(B\pi)$~\cite{BD}.
Also it was used to settle Gromov's conjecture in the almost spin case for virtual duality groups satisfying the coarse Baum-Connes conjecture~\cite{Dr1}.
There is no such technique available in the totally non-spin case, since neither the manifold nor its universal covering have a K-theory fundamental class.
This makes the totally non-spin case notoriously difficult. 

In this paper we prove Gromov's conjecture for virtual duality groups in the totally non-spin case with the exception when the dimension of a manifold equals the virtual cohomological dimension of the fundamental group, $\dim M\ne vcd(\pi_1(M))$. This is obtained as a corollary of the following theorem 
proven in the paper and the main result of~\cite{Dr1}.
\begin{thm}
 Let $M$ be a totally non-spin closed orientable $n$-manifold, $n\geq 5$,
 whose  universal cover $\Wi M$ satisfies the condition $\dim_{mc}\Wi M\le n-1$.
 Then $\dim_{mc}\Wi M\le n-2$.
\end{thm}

We recall that
every simply connected non-spin $n$-manifold, $n\ge 5$, admits a metric of positive scalar curvature. Perhaps one can conjecture that
every totally non-spin $n$-manifold $M$, $n\ge 5$, with (virtual) duality fundamental group $\pi$ admits a PSC metric whenever $n\ne vcd(\pi)$.
This is consistent with the main result of this paper which states that the inequality $\dim_{mc}\Wi M\le n-2$ conjectured by Gromov holds true for all  
manifolds $M$ with $\dim M\ne vcd(\pi_1(M))$.

We recall that virtual duality groups include large classes of groups such as
the virtually nilpotent groups, the arithmetic groups and the mapping class groups. 

It turns out that the Gromov conjecture in the case
$\dim M=vcd(\pi_1(M))$ is very special in PSC manifolds theory.  This equality holds for
aspherical manifolds where Gromov's conjecture implies the famous Gromov-Lawson conjecture:
\begin{conjec}[Gromov-Lawson]\label{G-L}
An aspherical manifold cannot carry a metric with positive scalar curvature.
\end{conjec}
The equality $\dim M=vcd(\pi_1(M))$ also presents in the following challenging problem in PSC theory which also
would be resolved by a proof of Gromov's conjecture.
\begin{question}\label{question}
 Does  the connected sum $M=T^{2n}\# \C P^n$ of the torus and complex projective space
 carry a metric with positive scalar curvature? 
\end{question}
It is reasonable to assume that the answer to this question is negative. 
In particular, the Gromov conjecture implies a negative answer to the question.
Since for odd $n$ the universal covering $\Wi M$ is spin and hypereuclidean,
$M$ does not admit a PSC metric by a theorem of Gromov and Lawson~\cite{GL}.
For even $n$ when $M$ is totally non-spin Question~\ref{question}  is a challenging problem. The minimal hypersurface method
of Schoen-Yau~\cite{SY} allows to treat the low dimensional cases when $n=2,4$.
There is a hope that their method can be extended to all dimensions~\cite{Lo}.

In this paper we reduce Gromov's conjecture for virtually abelian groups in the 
totally non-spin case to a version of the above problem. Also we note that
$T^{2n}\#\C P^n$ does not admit a metric of positive scalar curvature if the following conjecture
holds true for manifolds with abelian fundamental groups.
\begin{conjec}[$S^1$-Stability conjecture \cite{R3}]
A closed connected $n$-manifold $M$, $n>4$, admits a metric of positive scalar curvature if and only if $M\times S^1$ does.
\end{conjec}

The paper is arranged as follows. In \S 2 we present some facts about  PSC manifolds,
inessential manifolds and macroscopic dimension.
In \S 3 we prove the main result of this paper: The inequality  $$\dim_{mc}\Wi M\le n-2$$ for the universal covering of totally non-spin $n$-manifolds
whose fundamental group $\pi$ is a virtual FL duality group with $vcd(\pi)\ne n$.
In \S4 we investigate the case of (virtually) abelian fundamental groups.

\

In this paper we consider manifolds of dimension $\ge 5$. For 3-manifolds the Gromov conjecture was proved in~\cite{GL}. The case of 4-manifolds should be treated differently.

We are very grateful to Thomas Schick and Martin Nitsche who found flaws in the early versions of the paper.

\section{On  Inessential  Manifolds}

\subsection {Preliminaries}
Let $\pi=\pi_1(K)$ be the fundamental group of a CW complex $K$. By $u^K:K\to B\pi=K(\pi,1)$ we denote a map that classifies the universal covering $\Wi K$
of $K$. We refer to $u^K$
as a {\em classifying map} for $K$. We note that a map $f:K\to B\pi$ is a classifying map if and only if it induces an isomorphism of fundamental groups. 

The following theorem of Rosenberg is the main tool for dealing with Gromov's conjecture in the spin case.

\begin{thm}[Rosenberg~\cite{R1,R2}]\label{R}
Let $[M]_{KO}$ denote the fundamental class of a closed $spin$ $n$-manifold $M$  in the $KO$-theory. Let $\pi$ denote the fundamental group $\pi_1(M)$, then $\alpha\circ u^M_*([M]_{KO})=0$, where the homomorphism $u^M_* :KO_n(M) \to KO_n(B\pi)$ is induced by a classifying map $u^M:M\to B\pi$.
\end{thm}

The following result is the only known tool in the totally non-spin case.
\begin{thm}[Jung-Stolz~\cite{RS}]\label{JS}
Suppose that $N$ is a totally non-spin oriented manifold of dimension $\ge 5$ with $u^N_*([N])=u_*([M])$
in $H_n(B\pi)$  for some not necessarily 
connected oriented manifold  $M$ 
with a positive scalar curvature and a  map $u:M\to B\pi$.
Then $N$ admits a metric of positive scalar curvature.
\end{thm}

In this paper we use basic notations and facts from surgery theory and bordism theory~\cite{M},~\cite{Wa}.
We use the following
\begin{thm}[Surgery Theorem~\cite{GL},\cite{R3}]\label{surgeryTh}
Suppose that a manifold $N$ is obtained from a PSC manifold $M$ by a surgery of codimension $\le 3$.
Then $N$ admits a metric with positive scalar curvature.
\end{thm}

\subsection{Inessential manifolds and macroscopic dimension} We recall the following definition of Gromov~\cite{G3}:

\begin{defn}
An $n$-manifold $M$  with fundamental group $\pi$ is called {\it essential}  if its classifying map $u^M:M\to B\pi$ cannot be deformed into the $(n-1)$-skeleton $B\pi^{(n-1)}$ of a CW-complex structure on $B\pi$ and it is called {\em inessential} if $u^M$ can be deformed into $B\pi^{(n-1)}$.
\end{defn}

Note that for an inessential $n$-manifold $M$ we have $\dim_{mc}\Wi M\le n-1$. Indeed, a lift $\Wi{u^M}:\Wi M\to E\pi^{(n-1)}$ of a classifying
map is a uniformly cobounded proper map to an $(n-1)$-complex.
Generally, if a classifying map $u^M:M\to B\pi$ can be deformed to the $k$-dimensional skeleton, then $\dim_{mc}\Wi M\le k$.

Thus, one can consider a stronger version of Gromov's conjecture:
\begin{conjec}[The Strong Gromov Conjecture] A classifying map $u^M:M\to B\pi$ of the universal covering $\Wi M$ of  a closed  PSC  $n$-manifold $M$
with torsion free fundamental group can be deformed to the $(n-2)$-dimensional skeleton.
\end{conjec}
The restriction on the fundamental group is important,
since this conjecture is false for finite cyclic groups. For general groups one can consider
a virtual version of this conjecture. We note that in~\cite{BD} we proved the strong Gromov conjecture for products of free groups.

Thus, establishing the inessentiality of PSC manifolds  is the first step in a proof of the strong Gromov conjecture. 
We recall that the inessentiality of a manifold can be characterized as follows~\cite{Ba} (see also~\cite{BD}, Proposition 3.2). 
\begin{thm}\label{ref} Let $M$ be a closed oriented $n$-manifold. 
Then the following are equivalent:

1. $M$ is inessential;

2. $u^M_*([M])=0$ in $H_n(B\pi)$ where $[M]$ is the fundamental class of $[M]$.  
\end{thm}
In~\cite{BD} we proved the following addendum to Theorem~\ref{ref}.
\begin{prop}[\cite{BD}, Lemma 3.5]\label{ref2}
For an inessential manifold $M$  
with a CW complex structure a classifying map $u:M\to B\pi$ can be chosen such that 
$$u(M^{(n-1)})\subset B\pi^{(n-2)}.$$
\end{prop}

\subsection{Macroscopically inessential manifolds} 
The first step of Gromov's original conjecture is a statement about the macroscopic inessentiality of the universal covering of a closed PSC manifold.
One can split it off as a separate statement:

{\bf The Weak Gromov Conjecture.} {\it The macroscopic dimension of the universal covering $\Wi M$ of  a closed  PSC  $n$-manifold $M$ satisfies the inequality $$\dim_{mc}\Wi M\leq n-1$$ for a metric on $\Wi M$ lifted from $M$}.

\

The weak Gromov conjecture first appeared in~\cite{G1} in the language of filling radii.
Even the weak Gromov conjecture is out of reach,
since it implies the Gromov-Lawson conjecture (Conjecture~\ref{G-L}). 
The latter is known to be a Novikov type conjecture~\cite{R2}. 

\

There is an analog of Theorem~\ref{ref} for universal coverings.
\begin{thm}[\cite{Dr1}]\label{ref-mc}
Let $M$ be a closed oriented $n$-manifold  and let $\Wi u:\Wi M\to E\pi$ be a lift of $u^M$. 
Then the following are equivalent:

1. $\dim_{mc}\Wi M\le n-1$;

2. $\Wi u$ can be deformed by a bounded homotopy to $E\pi^{(n-1)}$;

3. $\Wi u_*([\Wi M])=0$ in $H_n^{lf}(E\pi)$ where $[\Wi M]$ is the fundamental class.
\end{thm}

REMARK. Formally, the condition 3 in Theorem~\ref{ref-mc} can be added only when $B\pi$ is special, in particular, when
$B\pi$ has finite $m$-skeletons. In fact only this case was used in ~\cite{Dr1} and in this paper. Generally, 
we abuse the notations and assume under the locally finite homology of $E\pi$
the direct limit $H^{lf}_n(E\pi)=\lim_{\rightarrow}H^{lf}_n(\Wi K_i)$ where $B\pi=\lim_{\rightarrow}K_i$
is the direct limit of finite complexes.   Since
the equivalence of conditions 1-2 holds true  for any choice of a CW complex $B\pi$, we prefer to state Theorem~\ref{ref-mc}
in full generality with this  modification of the definition of $H^{lf}_n(E\pi)$.

Thus the inequality  $\dim_{mc}\Wi M^n\le n-1$ is a macroscopic analog of inessentiality. We call manifolds
$N$ with $\dim_{mc}N<\dim N$  {\em macroscopically inessential}. Such manifolds are  called {\em macroscopically small} in~\cite{Dr1}.

There is an analog of Proposition~\ref{ref2}
\begin{prop}[\cite{Dr1}, Lemma 5.3.]\label{ref3}
For an $n$-dimensional manifold $M$ with a fixed CW complex structure, a classifying map $u:M\to B\pi$ and with macroscopically inessential universal covering $\Wi M$ 
any lift $\Wi u:\Wi M\to E\pi$ of $u$ admits a bounded deformation to a proper map $f:\Wi M\to E\pi^{(n-1)}$ with
$f(M^{(n-1)})\subset E\pi^{(n-2)}.$
\end{prop}
Here we assume that the metric on $E\pi$ is the lift of a proper metric on $B\pi$.

\subsection{Macroscopic dimension and QI-embeddings}
We recall that a map $f:X\to Y$ between metric spaces is
a quasi-isometric embedding (QI-embedding for short) if there are $\lambda,c>0$ such that the inequalities
$$
\frac{1}{\lambda}d_X(x_1,x_2)-c \le d_Y(f(x_1),f(x_2))\le \lambda d_X(x_1,x_2)+c.
$$ hold for all $x_1,x_2\in X$.
A map $f:X\to Y$ is a quasi-isometry if there is $D>0$ such that $f$ is a QI-embedding 
with a $D$-dense image.

\

Proposition~\ref{ref2} states in particular that for $n$-manifold $M$ the condition 
on the universal covering $\dim_{mc}\Wi M\le n-1$ implies that $\Wi M$ admits a continuous QI-embedding into $E\pi^{(n-1)}$
for $\pi=\pi_1(M)$. It turns out that the macroscopic dimension of universal coverings of closed manifolds can be characterized by
contiuous QI-embedding. To do a such characterization we make use of the following
\begin{prop}[\cite{Dr2}, Proposition 2.1]\label{md-charact}
For a proper  metric space $X$ the following are equivalent:

(1) $\dim_{mc}X\le n$;

(2) There is a continuous map $f:X\to K$ to a locally finite $n$-dimensional simplicial complex and a number $b>0$ such that 
$f^{-1}(\Delta)\ne\emptyset$ and $diam f^{-1}(\Delta)< b$ for every simplex $\Delta\subset K$.
\end{prop}
We recall that a metric space $E$ is called {\em uniformly $k$-connected} if for any $R>0$ there exists $S>0$ such that for each 
$i\le k$ the inclusion homomorphism of the $R$- and $S$-balls, $\pi_i(B_R(y))\to\pi_i(B_S(y))$, is trivial for all $y\in E$.
It is well-known that for a finite complex $B$ with $\pi_i(B)=0$ for $2\le i\le k$ its universal covering $E$ 
with a lifted metric is uniformly $k$-connected.
For a locally compact CW complex $B$ there is a relative version of the above
statement: 
\begin{prop}\label{unif}
Let $B$ be a CW complex given a proper metric with  finite $B^{(2)}$ and with $\pi_i(B)=0$ for $2\le i\le k$. Then for
any finite subcomplex $Y\subset B$, with $B^{(2)}\subset Y$ and a fixed $i\le k$ for any $R>0$  there exist a finite subcomplex $Z$,  $Y\subset Z\subset B$
and $S>0$ , such that  the inclusion homomorphism  $$\pi_i(B_R(y)\cap\Wi Y)\to\pi_i(B_S(y)\cap\Wi Z)$$ is trivial for all $y\in \Wi Y$
where the metric on $\Wi Z$ is lifted from $Z$.
\end{prop}
\begin{proof}
Let $p:\Wi Y\to Y$ be the universal covering map. We skip the cases when $i<2$ as obvious.
Note that the subgroup $G\subset \pi_i(Y)$ generated by $\cup_{y\in\Wi Y}p_*(\pi_i(B_R(y)))$ is a finitely generated as a $\pi_1(B)$-module.
Hence there is a compact $Z\subset B$ such that the inclusion homomorphism $\pi_i(Y)\to\pi_i(Z)$ takes $G$ to zero.
The maximal size of the $(i+1)$-disks killing the generators defines $S$.
\end{proof}

\begin{thm}\label{QI}
Suppose that  $E\pi$ is given a metric $d$ lifted from a proper metric on $B\pi$ where $\pi=\pi_1(X)$ for a finite CW complex $X$. 
Then for for its universal covering $\Wi X$ with the lifted metric  $\dim_{mc}\Wi X\le n$ if and only if there is
a continuous QI-embedding $g:\Wi X\to E\pi^{(n)}$.
\end{thm}
\begin{proof}
Let $f:\Wi X\to K$ be as in Proposition~\ref{md-charact} with $diam\, f^{-1}(\Delta)\le b$ for all simplices $\Delta$ in $K$.
Since $K^{(0)}\subset f(\Wi X)$,
we can define a map $\xi_0:K^{(0)}\to E\pi$ such that $\xi_0(v)=\Wi u(z)$ for some $z\in f^{-1}(v)$
where $\Wi u$ is a lift of a classifying map $u:X\to B\pi$. We may assume that $u(X)\subset B\pi^{(m)}$ with $m=\dim X$.
Note that $m\ge n$.

First, we assume that $B\pi^{(m)}$ is finite.
Then using uniform $(n-1)$-connectivity 
of $E\pi^{(m)}=\Wi{B\pi^{(m)}}$ with the lifted metric
by induction on $i$ we construct maps
$\xi_i:K^{(i)}\to E\pi^{(n)}$ extending $\xi_{i-1}$ with a uniform bound $diam\,\xi_i(\Delta)<D_i$ for all simplices $\Delta$ in $K^{(i)}$.
Let $\xi=\xi_n$. Thus,  $diam\,\xi(\Delta)<D=D_n$ for all simplices $\Delta$
in $K$.

We claim that the composition $g=\xi\circ f:\Wi X\to E\pi^{(n)}$ is a QI-embedding. Indeed, for $x_1, x_2\in\Wi X$ we take vertices $v_1\in \Delta_1$ and $v_2\in\Delta_2$ of simplices $\Delta_1$,
$\Delta_2$ containing $f(x_1)$ and $f(x_2)$ respectively. Then the triangle inequality and the fact that $\Wi u$ is a QI-embedding with some constants $\lambda$ and $c$
imply that
$$
d(g(x_1),g(x_2))\le d(\xi(v_1),\xi(v_2))+2D=d(\Wi u(z_1),\Wi u(z_2))+2D\le$$
$$\lambda d_{\Wi X}(z_1,z_2)+c+2D\le \lambda d_{\Wi X}(x_1,x_2)+2\lambda b+c +2D
$$
where $f(z_i)=v_i$ and, hence, $z_i,x_i\in f^{-1}(\Delta_i)$, $i=1,2$.
Similarly,
$$
d(g(x_1),g(x_2))\ge d(\xi(v_1),\xi(v_2))-2D=d(\Wi u(z_1),\Wi u(z_2))-2D\ge$$
$$
\frac{1}{\lambda} d_{\Wi X}(z_1,z_2)-c-2D\ge \frac{1}{\lambda} d_{\Wi X}(x_1,x_2)-2\frac{1}{\lambda}b-c -2D.
$$

If $B\pi^{(m)}$ is not finite we use Proposition~\ref{unif} to construct by induction a sequence of $\xi_i$ as above
together with a sequence of finite complexes
$Y_0\subset Y_1\subset\cdots Y_n\subset B$ such that $u(X)\subset Y_0$ and $p\xi_i(K^{(i)})\subset Y_i$
where $p:E\pi\to B\pi$ is the universal covering map.
The rest of the argument remains the same.

In the other direction, let $g:\Wi X\to E\pi^{(n)}$ be a QI-embedding with the constants $\lambda $ and $c$.
Note that for every closed $r$-ball $B_r(y)$ in $E\pi$ and any $x_1,x_2\in g^{-1}(B_r(y))$ we obtain
$$
\frac{1}{\lambda}d_{\Wi X}(x_1,x_2)-c\le d(g(x_1),g(x_2))\le 2r.
$$
Thus, $diam(g^{-1}(B_r(y))$ is closed and uniformly bounded. Hence $g$ is proper and uniformly cobounded.
\end{proof}
REMARK. In this paper we use only the easy part of Theorem~\ref{QI}.

\section{Gromov's conjecture for virtual duality groups}

We recall that the group of oriented relative bordisms    $\Omega_n(X,Y)$
of  the pair $(X,Y)$ consists of the equivalence classes of pairs $(M,f)$
where $M$ is an oriented $n$-manifold with boundary and $f:(M,\partial M)\to (X,Y)$ is continuous map.  Two pairs $(M,f)$ and $(N,g)$ are equivalent if there is a pair $(W,F)$, $F:W\to X$  called a {\em bordism} where $W$ is an orientable $(n+1)$-manifold with boundary such that  $\partial W =
M\cup W'\cup N$, $W'\cap M=\partial M$, $W'\cap N=\partial N$, $F|_M=f$,  $F|_N=g$, and $F(W')\subset Y$.  

In the special case when $X$ is one point, the manifold $W$ is called a bordism between $M$ and $N$.

\begin{prop}\label{bordism} For any CW complex $K$
there is an isomorphism $$\Omega_n(K,K^{(n-2)})\cong H_n(K,K^{(n-2)}).$$
\end{prop}
\begin{proof}
Since $\Omega_1(*)=0$ and
$K/K^{(n-2)}$ is $(n-2)$-connected, we obtain
that in the Atiyah-Hirzebruch spectral sequence on the diagonal  $p+q=n$ there is only one nonzero term 
which survives to $\infty$: $$E^2_{n,0}\cong E^{\infty}_{n,0}\cong
H_n(K,K^{(n-2)};\Omega_0(*))\cong
H_n(K,K^{(n-2)}).$$ Therefore,
$$\Omega_n(K,K^{(n-2)})\cong H_n(K,K^{(n-2)}).\eqno(*)
$$
\end{proof}

An $(n+1)$-dimensional $k$-handle is a space $H$ homeomorphic to the product $H\cong D^k\times D^{n+1-k}$. 
The subsets $D^k\times\{0\}\subset H$ and $\partial D^k\times D^{n+1-k}$  are called the core 
and  the base of the $k$-handle respectively. A $k$-handle $H$ is called attached to a $(n+1)$-manifold $W$ with boundary $\partial W$
if it intersects the boundary  along the base: $H\cap\partial W=\partial D^k\times D^{n+1-k}$ and $H\cap Int(W)=\emptyset$.  
A $k$-handle is also called a handle of index $k$.

We recall that for every bordism $W$ between $n$-manifolds $M$ and $N$ which is stationary on the boundary there is a handle decomposition $W=M\times[0,1]\cup\bigcup H_i\cup N\times[0,1]$ where $H_i\cong D^k\times D^{n+1-k}$. 
Moreover, for connected $M$ and $N$ there is a filtration 
$$M\times[0,1]=W_0\subset W_1\subset\dots \subset W_n=W\setminus (N\times(0,1])\subset W$$
where each $W_{i}$ is obtained from $W_{i-1}$ by attaching $i$-handles.
 Such a filtration defines a dual filtration
$$W\supset W\setminus (M\times[0,1))=W_n^*\supset\dots\supset W_1^*\supset W_0^*= N\times[0,1]$$ with the same set of handles
where each $i$-handle $H$ of $W_i$ is treated as a $(n-i+1)$-handle of $W_{n-i+1}^*$. 

This situation arises naturally for triangulated manifolds. Also it appears after a finite chain of surgeries.
In the paper we consider bordisms of open manifolds obtained from an infinite family of bordisms of compact manifolds with boundary
as follows. Let $M$ be an open $n$-manifold with a family of disjoint $n$-dimensional submanifold with boundary $\{V_\gamma\}$.
Let $U_\gamma$ be a family of stationary on the boundary bordisms between $V_\gamma$ and $N_\gamma$. Then the manifold $W=A\cup B$ with
$A=(M\setminus\cup_iInt\,V_\gamma)\times[0,1]$ and $B=\coprod_\gamma U_\gamma$ with $A\cap B=\coprod_\gamma\partial V_\gamma\times[0,1]$ is a {\em bordism 
defined by the family} $\{U_\gamma\}$. Thus $W$ is obtained from $M\times[0,1]$ by replacing the cylinders $V_\gamma\times[0,1]$ by bordisms $U_\gamma$.
Note that $\partial W=M\coprod N$. Thus, $W$ is a bordism
between manifolds $M$ and $N$. The $i$-handle filtrations on $U_\gamma$ define the filtration of $W$ 
$$M\times[0,1]=W_0\subset W_1\subset\dots \subset W_n=W\setminus (N\times(0,1])\subset W$$ and its dual.

For each $U_\gamma$ one can define a metric that extends the metric on $\partial V_\gamma\times[0,1]$ inherited from the
product metric on $M\times[0,1]$. This leads to a metric on $W$ such that the inclusion $M\subset W$ is an isometric embedding.
We call a bordism $W$ {\em bounded} if there is a uniform upper bound on the diameters of $U_\gamma$.

If $F:W\to X$ is a QI-embedding of a bounded bordism $W$ into a metric space $X$, then the pair $(W,F)$ is called a {\em bounded bordism} in $X$.
Note that in the case of a bounded $W$ the inclusion $M\to W$ is a quasi-isometry.

\begin{prop}\label{surgery1}
Let $W$ be a bounded bordism between an open $n$-manifolds $M$ and $N$  which is defined by a family of bordisms $\{V_\gamma\}$.
Suppose that $W$ does not have handles of dimension $\le k$.
Assume that $N$ admits a continuous QI-embedding $f:N\to K$ into a uniformly 
$(n-k-1)$-connected $l$-dimensional complex $K$. Then $\dim_{mc}M\le l$.
\end{prop}
\begin{proof}
Thus, there is a filtration
$$M\times[0,\epsilon]=W_0=\dots =W_k\subset W_{k+1}\subset\dots \subset W_n=W\setminus (N\times(0,\epsilon])\subset W$$
where each $W_{i}$ is obtained from $W_{i-1}$ by attaching $i$-handles  $H$ with $diam\, H< D$.
We consider the dual filtration
$$W\supset W\setminus (M\times[0,\epsilon))=W_n^*=\dots=W^*_{n-k}\supset\dots\supset W_1^*\supset W_0^*= N\times[0,\epsilon]$$

Since the projection $N\times[0,\epsilon]\to N$ is a quasi-isometry, we may assume that $f$ is defined on $N\times[0,\epsilon]$.
Since $K$ is uniformly $(n-k-1)$-connected and the family $f(N_\gamma\times[0,\epsilon])$ is uniformly bounded, the map $f$ can be extended to a map $g_1:W^*_1\to K$ with a uniform upper bound $R_1$ on the diameter of the images $g_1(H)$ of handles in $W_1^*$. This condition together with the assumption that the inclusion $W^*_0\subset W^*_1$ is a quasi-isometry imply that
$g_1$ is a QI-embedding.
Then $g_1$ can be extended to a quasi-isometric embedding $g_2:W^*_2\to K$ 
and so on. The uniform $(n-k-1)$-connectivity allows to proceed to a QI-embedding $g_{n-k}^*:W^*_{n-k}=W^*_n\to K$.
The projection $M\times[0,\epsilon]\to M$ defines a retraction $r:W\to W^*_n$. Since it is a quasi-isometry, the composition
$g=g^*_{n-k}\circ r$
is a QI-embedding. Since a continuous QI-embedding is a proper uniformly cobounded map (see Theorem~\ref{QI}), it follows that
$\dim_{mc}W\le l$ and, hence, $\dim_{mc}M\le l$.
\end{proof}

A map $f:X\to Y$ is called a $k$-equivalence if $f_*:\pi_i(X)\to\pi_i(Y)$ is an isomorphism for $i\le k$ and an epimorphism for $i=k+1$.
The proof of the following can be found in~\cite{Wa2}.
\begin{thm}\label{Wa}
Let $W$ be a bordism between compact manifolds $M$ and $N$ which is stationary on the boundary, $\partial M=\partial N$.
Suppose that the inclusion $M\to W$ is a $k$-equivalence. Then $W$ admits a handle decomposition with no handles of index $\le k+1$.
\end{thm}

Let  $\nu_X:X\to BSO$ denote
a classifying map for the stable normal bundle of  a  manifold $X$.

\begin{thm}\label{main3}  Let $M$ be a totally non-spin closed orientable $n$-manifold, $n\geq 5$,
with $\pi_1(M)=\Gamma$ whose  universal cover $\Wi M$ is  macroscopically inessential. Then a lift of  a classifying map $\Wi u^M:\Wi M\to E\Gamma$ can be boundedly deformed  to $E\Gamma^{(n-2)}$, in particular, $\dim_{mc}\Wi M\le n-2$. 
\end{thm}
\begin {proof}
The assumption that $\Wi M$ is non-spin implies that  $$(\nu_{\Wi M})_*:\pi_2(\Wi M)\to\pi_2(BSO)=\Z_2$$ is surjective.
Let $\Sigma'$ be a 2-sphere in $\Wi M$ with $(\nu_{\Wi M})_*([\Sigma'])\ne 0$ where $[\Sigma]$ denotes the corresponding element of $\pi_2(\Wi M)$.
Since the covering map $\Wi M\to M$ induces an isomorphism of 2-dimensional homotopy groups and $dim\, M\ge 5$, 
we may assume that $\Sigma'$ is a lift of an embedded sphere $\Sigma\subset M$.
Let $V$ be a regular neighborhood of $\Sigma$ with a lift to $V'$, a regular neighborhood of $\Sigma'$.
Let $V_{\gamma}=\gamma(V')$ denote a $\gamma$-translate of $V'$.
In view of the uniqueness of lifting of $V$ with a prescribed lift of a point, the family $V_{\gamma}$ is  disjoint.

We assume that the CW structure on $M$ has one $n$-dimensional cell. In view of Proposition~\ref{ref3}  there is a bounded deformation of $\Wi u$ to a map $f:\Wi M\to E\Gamma^{(n-1)}$
with $f(\Wi M\setminus \coprod_{\gamma\in\Gamma}D_{\gamma})\subset E\Gamma^{(n-2)}$ where $\{D_{\gamma}\}_{\gamma\in\Gamma}$ are the lifts of a fixed closed $n$-ball $D$ in 
the top dimensional cell in $M$. 
Clearly, $f$ is a QI-embedding.
We may assume that $D\subset V$.

Note that the restriction
of $f$ to $(D_{\gamma},\partial D_{\gamma})$ defines a zero element in $H_n(E\Gamma,E\Gamma^{(n-2)})$. Moreover, there is $r>0$ such that $f(D_{\gamma})\subset B_r(f(c_{\gamma}))$ where $c_{\gamma}\in D_{\gamma}$ and
$f|_{D_{\gamma}}$ defines a zero element in $$H_n(B_r(f(c_{\gamma})), B_r(f(c_{\gamma}))\cap E\Gamma^{(n-2)}).$$

By Proposition~\ref{bordism} there is a relative bordism
$(W_{\gamma},q_{\gamma})$ of $(D_{\gamma},\partial D_{\gamma})$ to $(N_{\gamma},S'_{\gamma})$ with $q_{\gamma}(N_{\gamma}\cup\partial W_{\gamma}\setminus D_{\gamma})\subset E\Gamma^{(n-2)}$.
We may assume that the bordism $W'\subset\partial W_{\gamma}$ of the boundaries $\partial D_{\gamma}\cong S^{n-1}$ and $S'_{\gamma}$ is stationary, $W'\cong\partial D_{\gamma}\times[0,1]$
and $q(x,t)=q(x)$ for all $x\in\partial D$ and all $t\in[0,1]$. 
By performing 1-surgery on $W_{\gamma}$ we may assume that $W_{\gamma}$ is simply connected.

For each $\gamma\in\Gamma$ we can enlarge the  bordism $W_{\gamma}$ 
by the trivial bordism $(V_{\gamma}\setminus Int(D_{\gamma}))\times[0,1]$ to a bordism
$\bar W_{\gamma}$ of manifolds with boundaries $V_{\gamma}$ and $N_{\gamma}'$. 
The map $q_{\gamma}$  extends to a map $\bar q_{\gamma}:\bar W_{\gamma}\to E\Gamma$ with $\bar q_{\gamma}(V_{\gamma}\setminus Int(D_{\gamma}))\times[0,1])\subset E\Gamma^{(n-2)}$ by means of $f$.
Thus $(\bar W_{\gamma},\bar q_{\gamma})$ is a stationary on the boundary bordism between 
$V_{\gamma}$ and $N'_{\gamma}$ with $\bar q_{\gamma}(N'_{\gamma})\subset E\Gamma^{(n-2)}$.

Note that the inclusion $V_{\gamma}\to\bar W_{\gamma}$
induces an isomorphism of the fundamental groups. By Theorem~\ref{Wa} we may assume that $\bar W_{\gamma}$ does not have handles in dimension 1.

Note that every 2-sphere $S$ that generates an element of the kernel
of $(\nu_{\bar W_{\gamma}})_*:\pi_2(\bar W_{\gamma})\to\pi_2(BSO)$, has trivial stable normal bundle.
Hence we can apply
a  surgery in dimension 2 on $\bar W_{\gamma}$ to obtain a manifold $\hat W_{\gamma}$ and a map $\nu_{\hat W_{\gamma}}:
\hat W_{\gamma} \to BSO$ that induces a monomorphism $(\nu_{\hat W_{\gamma}})_*:\pi_2(\hat W_{\gamma})\to\pi_2(BSO)=\Z_2$.
The maps $\bar q_{\gamma}$ can be modified to a maps $\hat q_{\gamma}:\hat W_{\gamma}\to E\Gamma$ with  $\hat q_{\gamma}=\bar q_{\gamma}$
on $\partial \hat W_{\gamma}=\partial\bar W_{\gamma}$ and with a uniform bound for the distance between the images $im(\hat q_{\gamma})$ and $im(\bar q_{\gamma})$. Let $i^\gamma:V_{\gamma}\to \hat W_{\gamma}$ denote the inclusion.
Since $(\nu_{V_{\gamma}})_*$ is surjective, $(\nu_{\hat W_{\gamma}})_*\circ i_*^\gamma=\nu_{V_{\gamma}}$
and $(\nu_{\hat W_{\gamma}})_*$ is injective, we obtain that $i_*^\gamma:\pi_2(V_{\gamma})\to
\pi_2(\hat W_{\gamma})$ is surjective.
Hence,  by Theorem~\ref{Wa} we may assume that $\bar W_{\gamma}$ does not have handles in dimension 2.

Let $\hat W$ be a bordism of $\Wi M$ defined by the family of relative bordisms $\coprod_{\gamma\in\Gamma}\hat W_{\gamma}$ with $\partial W=\Wi M\coprod N$.  Let $i: \Wi M\to\hat W$ denote the inclusion map.
We may choose a metric on $\hat W$ such that $\hat W$ is bounded and $i$ is an isometric embedding. 
The union of the maps $\hat q_{\gamma}$ naturally extends to a QI-embedding $\hat q:\hat W\to E\Gamma$.

We note that the bordism $(\hat W,\hat q)$ with the continuous QI-embedding $\hat q:\hat W\to E\Gamma$  between $(\Wi M,f)$ and $(N,g)$
is a bounded bordism  which does not have handles of dimension $\le 2$ with a QI-embedding
$g:N\to E\Gamma^{(n-2)}$. In the case when $B\pi^{(n-1)}$ is finite the space $E\Gamma^{(n-2)}$
is uniformly $(n-3)$-connected. Then Proposition~\ref{surgery1} implies that $\dim_{mc}\Wi M\le n-2$.
In the general case we use induction and Proposition~\ref{unif} to construct subcomplexes $Y_1\subset\dots\subset Y_{n-2}\subset B\pi$ and QI-embeddings $g_i:W^*_i\to \Wi Y_i^{(n-2)}$ to complete the proof.
\end{proof}

We recall that the  groups that admit finite $B\pi$ are called {\em geometrically finite}
or of the type $FL$.  A group $\pi$ is of the type $FP$ if $B\pi$ is dominated by a finite complex.
It is  an open problem whether $FP=FL$~\cite{Br}.

We recall that a group $\pi$  is called a {\em duality group}~\cite{Br} 
if it is of the type $FP$ and there is a $\pi$-module $D$ such that
$$
H^i(\pi,M)\cong H_{m-i}(\pi,M\otimes D)
$$
for all $\pi$-modules $M$ and all $i$ where $m=cd(\pi)$ is the cohomological dimension of $\pi$.

A group $\pi$ is a virtual FL duality group if it contains a finite index subgroup $\pi'$ which is a FL duality group.

\begin{thm}\label{main}
Gromov's conjecture holds true for manifolds $M$ whose fundamental groups $\pi=\pi_1(M)$ are virtual FL duality groups if one of the following holds

(1) $M$ is almost spin and $\pi$ satisfies the coarse Baum-Connes conjecture;

(2) $M$ is totally non-spin and $\dim M\ne vcd\,\pi$.
\end{thm}
\begin{proof}
(1) is Theorem 5.6 from~\cite{Dr1}.

(2) We may assume that $M$ is orientable. Othrewise we consider an orientable 2-to-1 covering $M'\to M$ and note that they have the same universal cover $\Wi M'=\Wi M$.
Let $\pi'$ be a finite index subgroup of $\pi$ that is a FL duality group
and let $M'\to M$ be corresponding covering. Then there is a finite CW complex $B\pi'$. In that case
by Proposition 5.5~\cite{Dr1}, $H_n^{lf}(E\pi';\Z)=0$
for $n\ne cd(\pi')=vcd(\pi)$.  Theorem~\ref{ref-mc} implies that the universal cover $\Wi M'=\Wi M$  is macroscopically inessential provided $\dim M\ne vcd(\pi)$.
Then by Theorem~\ref{main3} $\dim_{mc}\Wi M\le n-2$.
\end{proof}

\section{The case of abelian fundamental group}

Theorem ~\ref{main} implies the following fact which also can be derived from~\cite{Dr1}.

\begin{thm}
Suppose that a closed oriented $n$-manifold $M$ with positive scalar curvature has the virtually free abelian fundamental group
$\pi_1(M)$ of rank $r\ne n$. Then $\dim_{mc}\Wi M<n$.
\end{thm}
In particular, the weak Gromov conjecture holds true under the condition $r\ne n$. The case
$r=n$ should be treated differently: If the manifold is almost spin
then $\dim_{mc}\Wi M<n-1$ by \cite{Dr1}, Theorem 5.6, since virtually abelian groups satisfy the coarse Baum-Connes conjecture.
For the totally non-spin case we have only a reduction of the conjecture to a version of
Question~\ref{question} and the $S^1$-stability conjecture.

\subsection{$S^1$-stability conjecture} We recall that $H_m(X)^+$ denotes the subset of integral homology classes which can be realized by orientable manifolds with positive scalar curvature. It is known that $H_m(X)^+\subset H_m(X)$ is a subgroup~\cite{R3}.

\begin{thm}[T. Schick~\cite{Sch}]
Let $\alpha\in H^1(X)$. Then  for $3\le k\le 8$ the cap product with $\alpha$ takes $H_k(X)^+$ to $H_{k-1}(X)^+$.
\end{thm}
\begin{prop}\label{S-1}
The $S^1$-stability conjecture for manifolds with free abelian fundamental group implies that $H_*(T^n)^+=0$, $n\geq 5$.
\end{prop}
\begin{proof} Let $p:\#_rT^k\to T^k$ be a map of degree $r$.
Denote by $rT^k_s$ the $k$-manifold obtained by 1-surgery from the connected sum $\#_rT^k$ of $r$ copies of $T^k$ that kills the kernel
of the homomorphism $$p_*:\pi_1(\#_rT^k)\to\pi_1(T^k).$$
The surgery changes the map $p$ into  a classifying map $u:rT^k_s\to T^k$  with $deg(u)=r$.
Suppose that $r[T^n]\in H_*(T^n)^+$, $n\ge 5$. By the Jung-Stolz theorem $(rT^4_s\#\C P^2)\times T^{n-4}$ has a positive scalar curvature metric. Then by the $S^1$-stability conjecture, $M=rT^4_s\# \C P^2$
has such a metric. We apply Schick's theorem  two times consecutively with the  1-dimensional cohomology  classes $\bar\alpha_1$ and $\bar\alpha_2$ generated by
the collapse $u:M\to T^4$ of $M$ onto $T^4$ followed by the projections onto the first factor $S^1$  and the second factor respectively. 
This produces a surface $j:S\to M$
with positive scalar curvature which realizes the 2-homology class $[M]\cap(\bar\alpha_1\cup\bar\alpha_2)$. Thus, $S$ is the 2-sphere.
On the other hand, $u_*(j_*[S])=[T^4]\cap(\alpha_1\cup\alpha_2)\ne 0$ where $\alpha_1$ and $\alpha_2$ are generators of $H^1(T^4)$ generated by projections to the first and second factors. This produces a
contradiction since every mapping of the 2-sphere into a torus is nullhomotopic.
\end{proof}
\begin{question}
Is the condition $H_n(T^n)^+\ne 0$ equivalent to the equality  $H_n(T^n)^+=H_n(T^n)$ ?
\end{question}
\begin{prop}\label{+}
For each $n$ the following conditions are equivalent

(1) $H_{4n}(T^{4n})^+=H_{4n}(T^{4n})$;

(2) $T^{4n}\#\C P^{2n}$ admits a positive scalar curvature metric.
\end{prop}
\begin{proof}
(1) $\Rightarrow$ (2) follows from the Jung-Stolz theorem since $T^{4n}\#\C P^{2n}$ is totally non-spin.

(1) $\Leftarrow$ (2) follows from the fact that $H_*(X)^+$ is a group.
\end{proof}

\subsection{Connection to inessentiality}
Since the $n$-dimensional homology the $m$-dimensional torus is generated by $n$-dimensional subtori,
we obtain the following:
\begin{prop}\label{projection}
Suppose that a map $f:M\to T^m$ of a closed oriented $n$-manifolds takes the fundamental class to
a nonzero element in $H_n(T^m)$. Then there is a projection $q:T^m\to T^n$ such that $deg(q\circ f)\ne 0$.
\end{prop}
\begin{proof} There exists an integral $n$-dimensional cohomology class $\beta$ with nonzero evaluation 
$\langle f_*([M]),\beta\rangle$.
Recall that the cohomology ring of the torus $T^m=S^1\times\dots\times S^1$ is the exterior algebra
$H^*(T^m)=\Lambda[\bar\alpha_1,\dots,\bar\alpha_m]$ on 1-dimensional generators that come from 
the generators $\alpha_i$ of the factors. Thus, 
$$\beta=\sum_{I=(i_1,\dots, i_n)}c_I\bar\alpha_{i_1}\wedge\dots\wedge\bar\alpha_n.$$
Therefore $a=\langle f_*([M]),\bar\alpha_{j_1}\wedge\dots\wedge\bar\alpha_{j_n}\rangle\ne 0$ for some
$j_1,\dots, j_n$. We define $q:T^m\to T^n$ to be the projection onto the factor $S^1_{j_1}\times\cdots\times S^1_{j_n}$.
Then $$ a=
\langle f_*([M]),q^*(\alpha_{j_1}\wedge\dots\wedge\alpha_{j_n})\rangle=\langle q_*f_*([M]),\alpha_{j_1}\wedge\dots\wedge\alpha_{j_n}\rangle\ne 0.$$ Hence $q_*f_*([M])\ne 0$.
\end{proof}

\begin{thm}\label{main1}
Let $M$ be closed orientable $k$-manifold that admits a metric of positive scalar curvature with $\pi_1(M)=\Z^m$.
Suppose that $H^{k}(T^{k})^+=0$. Then  $M$ is  inessential.
\end{thm}
\begin{proof}
Let $u^M:M\to T^m$ be a classifying map. Assume that $u^M_*([M])\ne 0$. By Proposition~\ref{projection}
there is $q:T^m\to T^k$ such that $(q\circ u^M)_*([M])\ne 0$. Then  $H_k(T^k)^+\ne 0$. We obtain a contradiction.
\end{proof}

\subsection{Deformation into $B\pi^{(n-2)}$}

Proposition 4.8 and Theorem 4.9 below correspond to Proposition~\ref{surgery1} and Theorem~\ref{main3} respectively.

\begin{prop}\label{surgery}
Suppose that a closed $n$-manifold $N$ is obtained from a manifold $M$ by a chain of $k$-surgeries with $2\le k\le n-1$
such that the inclusion of $N$ to the corresponding bordism $W$ induces an isomorphism of the fundamental groups.
Assume that a classifying map $u^N:N\to B\pi$ admits a deformation to $B\pi^{(n-2)}$. Then so does $u^M:M\to B\pi$.
\end{prop}
\begin{proof}
The conditions on the bordism $W$ imply that  $W$ admits the following dual decomposition
$$W\supset W\setminus (M\times[0,1))=W_{n-2}^*\supset\dots\supset W_3^*\supset W_2^*\supset W^*_1= N\times[0,1]$$ 
where $W_2^*$ is obtained from $W_1^*$ by attaching 2-cells by means of
null-homotopic maps and thickening, $W_3^*$ is obtained from $W_2^*$ by attaching 3-cells and thickening and so on up to
attaching $(n-2)$-dimensional cells and thickening.  
Therefore,  the map $u^N$ can be extended to a map $g:W\to B\pi^{(n-2)}$.
Since the inclusion $M\to W$ is a 1-equivalence we obtain that the restriction $g|_{M}$ induces isomorphism of the fundamental groups, and  hence is a classifying map for $\Wi M$.
\end{proof}

\begin{thm}\label{main2} Let $M$ be a totally non-spin closed orientable inessential $n$-manifold, $n\geq 5$. Then a classifying map $u^M:M\to B\pi$ can be deformed  to $B\pi^{(n-2)}$, in particular, $\dim_{mc}\Wi M\le n-2$. 
\end{thm}
\begin{proof}
Here we use an accordingly modified argument of Theorem~\ref{main3}.

We assume that the CW structure on $M$ has one $n$-dimensional cell. In view of Lemma~\ref{ref2}  there is a classifying map $f:M\to B\pi$
with $f(M\setminus D)\subset B\pi^{(n-2)}$ where $D$ is an $n$-ball.
Note that the restriction
of $f$ to $(D,\partial D)$ defines a zero element in $H_n(B\pi,B\pi^{(n-2)})$. Therefore by Proposition~\ref{bordism} there is a relative bordism
$(W,q)$ of $(D,\partial D)$ to $(N',S)$ with $q(N'\cup\partial W\setminus D)\subset B\pi^{(n-2)}$.
We may assume that the bordism $W'\subset\partial W$ of the boundaries $\partial D\cong S^{n-1}$ and $S$ is stationary, $W'\cong\partial D\times[0,1]$
and $q(x,t)=q(x)$ for all $x\in\partial D$ and all $t\in[0,1]$. 
Let $\bar W$ be an extension of $W$ to a bordism of $M$ by the product bordism.  Let $\bar q: \bar W\to B\pi$ denote the extension of $q$
by means of $f|_{M\setminus D}$. Thus $\bar W$ is a bordism between $M$ and $N$. 

We may assume that the restriction of $\bar q$ to $N$
induces an isomorphism of the fundamental groups. If not, we can modify  $\bar W$  by attaching 1-handles and 2-handles to $N$
to achieve this.
By performing 1-surgery on $\bar W$ that kills the kernel of $\bar q_*:\pi_1(W)\to\pi_1(B\pi)=\pi$
we may assume that $\bar q$ induces an isomorphism of the fundamental groups.
Let $i_M: M\to \bar W$ and $i_N:\to\bar W$ denote the inclusion maps.
Thus, both $i_M$ and $i_N$ induce  isomorphisms of the fundamental groups. 

Note that every 2-sphere $S$ that generates an element of the kernel
of $(\nu_{\bar W})_*:\pi_2(\bar W)\to\pi_2(BSO)$ has trivial stable normal bundle.
It is easy to see that $\pi_2(\bar W)$ as a $\pi$-module is finitely generated. Hence we can apply
surgery in dimension 2 on $\bar W$ to obtain a map $\nu_{\bar W}:
\bar W \to BSO$ that induces an isomorphism $(\nu_{\bar W})_*:\pi_2(\bar W)\to\pi_2(BSO)$.
The assumption that $\Wi M$ is non-spin implies that  $(\nu_{M})_*:\pi_2(M)\to\pi_2(BSO)=\Z_2$ is surjective.
Since $\nu_M=\nu_{\bar W}\circ i_M$ is an epimorphism and $(\nu_{\bar W})_*$ is an isomorphism, it follows that 
$i_*:\pi_2(M)\to\pi_2(\bar W)$ is  an epimorphism.
Therefore (see Theorem~\ref{Wa}), $W$ is a bordism between $M$ and $N$ which is obtained from $M$ as the result  of $k$-surgeries with $2\le k\le n-1$.
Note that $g=\bar q|_N:N\to B\pi^{(n-2)}$ induces an isomorphism of the fundamental groups and, hence, is a classifying map.
 Proposition~\ref{surgery} completes the proof.
\end{proof}

\begin{thm}\label{virtfree}  The strong
Gromov conjecture holds true for totally non-spin oriented $n$-manifolds $M$, $n\ge 5$, with free abelian fundamental group if and only if $H_n(T^n)^+=0$ . 
\end{thm}
\begin{proof} Suppose that $H_n(T^n)^+=0$ .
Then by Theorem~\ref{main1} $M$ is inessential.
Theorem~\ref{main2} implies that $u^M$ can be deformed to the $(n-2)$-skeleton.

Let $f:M\to T^n$ be a  map  of a PSC manifold with $f_*([M])\ne 0$. We perform 0 and 1 surgery on $M$ to obtain
a manifold $N$  with a classifying map $u:N\to T^n$ such that $u_*([N])=f_*([M])$. By Theorem~\ref{surgeryTh} $N$ admits a PSC metric.
By the strong Gromov conjecture $u$ should be deformable to the $(n-2)$-dimensional skeleton. In particular, $N$s inessential.
This contradicts to Theorem~\ref{ref}.
\end{proof}
By going to  finite coverings, we derive the following
\begin{cor}
Suppose that $H_n(T^n)^+=0$. Then 
Gromov's conjecture holds true for $n$-manifolds $M$, $n\ge 5$ with virtually free abelian fundamental group. 
\end{cor}
REMARK. The case of non-orientable $M$ can be reduced to the orientable case by a 2-to-1 covering.

\begin{cor}
The $S^1$-stability conjecture implies the Gromov conjecture for virtually abelian groups for orientable manifolds of dimension $\geq 5$.
\end{cor}

\end{document}